\documentclass[twoside]{amsart}
\usepackage{amsmath}
\usepackage{amsfonts}
\usepackage{amssymb,enumerate}
\usepackage{amsthm}
\usepackage{hyperref}

\theoremstyle{plain}
\newtheorem{theorem}{Theorem}[section]
\newtheorem{lemma}[theorem]{Lemma}
\newtheorem{definition}[theorem]{Definition}

\begin{document}

\title{  Hierarchies of Subsystems of Weak Arithmetic}

\author{ Shahram Mohsenipour}
\address{\ Shahram Mohsenipour,
         School of Mathematics, Institute for Research in Fundamental Sciences (IPM)\\
         P. O. Box 19395-5746, Tehran, Iran}
\email{mohseni@ipm.ir}
\thanks{This work was done while the author was a Postdoctoral Research Associate at the School of Mathematics,
        Institute for Research in Fundamental Sciences (IPM)}


\begin{abstract}
We completely characterize the logical hierarchy of
various subsystems of weak arithmetic, namely: ZR, ZR + N, ZR +
GCD, ZR + Bez, OI + N, OI + GCD, OI + Bez.
\end{abstract}

\maketitle

\bibliographystyle{amsplain}

\section{Introduction}

In 1964 Shepherdson \cite{shep} introduced a weak system of
arithmetic, Open Induction (OI), in which the Tennenbaum phenomenon
does not hold. More precisely, if we restrict induction just to open
formulas (with parameters), then we have a recursive nonstandard
model. Since then several authors have studied Open Induction and
its related fragments of arithmetic. For instance, since Open
Induction is too weak to prove many true statements of number theory
(It cannot even prove the irrationality of $\surd{2}$), a number of
algebraic first order properties have been suggested to be added to
OI in order to obtain closer systems to number theory. These
properties include: Normality \cite{van} (abbreviated by N), having
the GCD property \cite{smith}, being a Bezout domain \cite{MM,smith}
(abbreviated by Bez), and so on. We mention that GCD is stronger
than N, Bez is stronger than GCD and Bez is weaker than $IE_{1}$
($IE_{1}$ is the fragment of arithmetic based on the induction
scheme for bounded existential formulas and by a result of Wilmers
\cite{wilm}, does not have a recursive nonstandard model).
Boughattas in \cite{b,Bb} studied the non-finite axiomatizability problem
and established several new results, including: (1) OI is not
finitely axiomatizable, (2) OI + N is not finitely axiomatizable. To
show that, he defined and considered the subsystems (OI)$_{p}$ of
(OI) and (N)$_{n}$ of N ($1\leq p,n <\omega$) (See the next section
for the definitions) and proved:

\begin{theorem}[Boughattas \cite{b}]

$(1)$ $(OI)_{p}$ is finitely axiomatizable,

$(2)$ Suppose $(p!,p^{'})=1$, then $(OI)_{p}\nvdash(OI)_{p^{'}}$.
\end{theorem}

\begin{theorem}[Boughattas \cite{Bb}, Theorem 2]Suppose
$(p!,p^{'})=1$ and $(n!,n^{'})=1$,

$(1)$ $N+(OI)_{p}\nvdash(OI)_{p^{'}}$,

$(2)$ $(N)_{n}+(OI)\nvdash (N)_{n^{'}}$,

$(3)$ $(OI)_{p}+\neg(OI)_{p^{'}}+(N)_{n}+ \neg (N)_{n^{'}}$ is
consistent.
\end{theorem}

In \cite{moh2} we strengthened Theorem 1.1 (2) to completely
characterize the logical hierarchy of OI, by showing that
$(OI)_{p}\nvdash(OI)_{p+1}$ iff $p\neq3$. In this paper by modifying
Boughattas' original proofs, we also strengthen Theorem 1.2 in two
directions and completely characterize
the logical hierarchy of OI + N, OI + GCD, OI + Bez:\\

\noindent{\bf Theorem C.} $Bez$ + $(OI)_{p}$ $\nvdash$
$(OI)_{p+1}$, {\em when} $p\neq 3$.\\

\noindent{\bf Theorem D.} $(OI)_{p}$ + $\neg(OI)_{p+1}$ +
$(N)_{n}$ + $\neg(N)_{n+1}$ {\em is consistent, when} $p\neq 3$.\\

So we will have the following immediate consequences:\\

\noindent{\bf Corollary E.}

(1) $N$ + $(OI)_{p}$ $\nvdash$ $(OI)_{p+1}$, {\em when} $p\neq 3$.

(2) $GCD$ + $(OI)_{p}$ $\nvdash$ $(OI)_{p+1}$, {\em when} $p\neq
3$.

(3) {\em All of the following subsystems of arithmetic are
non-finite axiomatizable: OI, OI + N, OI + GCD,
OI + Bez, (OI)$_{p}$ + N, OI + (N)$_{n}$}.\\

In Theorems A and B of this paper, we consider the ZR versions of
the above theorems. ZR is a subsystem of arithmetic that allows
Euclidean division over each non-zero natural number
$n\in\mathbb{N}$. ZR is introduced by Wilkie \cite{wil} in which he
proved that ZR and OI have the same $\forall_{1}$-consequences.
Later developments showed that ZR had very important role in
constructing models of OI (See Macintyre-Marker \cite{MM}, Smith
\cite {smith}). ZR + N has also been studied in \cite{o}. In Theorem
A, we study natural subsystems (ZR)$_{S}$ of (ZR), for a nonempty
subset $S$ of the set of prime numbers $\mathbb{P}$ (see the next
section for
definition) and show that:\\

\noindent {\bf Theorem A.} {\em Suppose S is a nonempty subset of
$\mathbb{P}$ and $q$ is a prime number such that $q\notin S$, then
$(ZR)_{S}$ + $Bez$ $\nvdash$
 $(ZR)_{q}$.}\\

Boughattas in (\cite{Bb}, Lemma 5) proved that DOR + N and ZR + N
are not finitely axiomatizable. More precisely he showed that:

\begin{theorem}[Boughattas \cite{Bb}, Lemma 5] Suppose
$(n!,n^{'})=1$. Then $ZR+(N)_{n}\nvdash(N)_{n^{'}}$.
\end{theorem}

We modify Boughattas' proof and strengthen the above theorem in
Theorem B:\\

\noindent{\bf Theorem B.} {\em Suppose S is a nonempty subset of
$\mathbb{P}$ and $q$ is a prime number such that $q\notin S$, then
$(ZR)_{S} + (N)_{n} + \neg (ZR)_{q} + \neg (N)_{n+1}$ is
consistent.}\\

Therefore we will have the following immediate implications:\\

\noindent{\bf Corollary F.} {\em Suppose S is a nonempty subset of
$\mathbb{P}$ and $q$ is prime number such that $q\notin S$, then}

(1) $(ZR)_{S}$ + $N$ $\nvdash$
 $(ZR)_{q}$.

(2) $(ZR)_{S}$ + $GCD$ $\nvdash$
 $(ZR)_{q}$.

(3) {\em All of the following subsystems of arithmetic are
non-finite axiomatizable: ZR, ZR + N, ZR + GCD, ZR + Bez, (ZR)$_{S}$
+ N, ZR + (N)$_{n}$, when S is an infinite subset of the set of
prime numbers}.

\section{Preliminaries}
Let $L$ be the language of ordered rings based on the symbols +,
$-$, $\cdot$, 0, 1, $\leq$. We write $\mathbb{N^\ast}$ for
$\mathbb{N}\setminus\{0\}$. We will work with the following set of
axioms in $L$:\\

$\bf DOR$: discretely ordered rings, i.e., axioms for ordered
rings and
\begin{center}
$\forall x \neg(0<x<1)$.
\end{center}

$\bf ZR$: discretely ordered $\mathbb{Z}$-rings, i.e., DOR and for
every $n\in\mathbb{N^{\ast}}$
\begin{center} $\forall x\exists q,r (x=nq+r\bigwedge 0\leqslant r<n)$.
\end{center}

We denote the sentence ``DOR + $\forall x\exists q,r
(x=nq+r\bigwedge 0\leqslant r<n)$" by (ZR)$_{n}$. Suppose
$\mathbb{P}$ denote the set of prime numbers of $\mathbb{N}$. Let
$S$ be a nonempty subset of $\mathbb{P}$. We
define the subsystem (ZR)$_{S}$ of ZR as the below:\\

${\bf ZR}_{S}$:  DOR + for every $p\in S$
\begin{center} $\forall x\exists q,r (x=qp+r\bigwedge 0\leqslant r<p)$.
\end{center}

If $S=\{p_{i_{1}},\ldots,p_{i_{n}}\}$ is a finite subset of
$\mathbb{P}$, we write (ZR)$_{p_{i_{1}},\ldots,p_{i_{n}}}$ instead
of (ZR)$_{\{p_{i_{1}},\ldots,p_{i_{n}}\}}$. This is consistent
with the above
notation (ZR)$_{n}$.\\

$\bf {OI}$: open induction, i.e., DOR and for every open
$L$-formula $\psi(\bar{x},y)$
\begin{center}
$\forall \bar{x} (\psi(\bar{x},0) \bigwedge \forall y\geqslant 0
(\psi(\bar{x},y)\rightarrow\psi(\bar{x},y+1)) \rightarrow\forall
y\geqslant 0  \psi(\bar{x},y))$.\\
\end{center}

By considering the fact that in discretely ordered rings an open
$L$-formula $\varphi(\bar{x},y)$ can be written as a Boolean
combination of polynomial equalities and inequalities with the
variable $y$ and the parameters $\bar{x}$, there exist natural
numbers $m,n$ such that:

\begin{center}
$\varphi(\bar{x},y)=\bigwedge_{i\leq m }\bigvee_{j\leq n
}p_{ij}(\bar{x},y)\leq q_{ij}(\bar{x},y),$
\end{center}
we can define the degree of $\varphi(\bar{x},y)$ relative to $y$
by

\begin{center}
deg $\varphi(\bar{x},y)$= max $\{$deg$_{y}$ $p_{ij}(\bar{x},y)$,
deg$_{y}$
$q_{ij}(\bar{x},y)|i\leq m, j\leq n\}.$\\
\end{center}

$({\bf OI})_p$ : open induction up to degree $p$ (i.e., DOR and
for every open $L$-formula $\psi(\bar{x},y)$ with deg
$\psi(\bar{x},y)\leq p$

\begin{center}
$\forall \bar{x} (\psi(\bar{x},0) \bigwedge \forall y\geqslant 0
(\psi(\bar{x},y)\rightarrow\psi(\bar{x},y+1)) \rightarrow\forall
y\geqslant 0  \psi(\bar{x},y)))$.\\
\end{center}

$\bf N$: normality (i.e., being domain and integrally closed in
its fraction field, namely for every $n\in\mathbb{N}^{*}$,
$\forall x, y, z_1, \ldots, z_{n}$

\hspace{.8in}$(y\neq 0\wedge
x^{n}+z_{1}x^{n-1}y+\ldots+z_{n-1}xy^{n-1}+z_{n}y^{n}=0\longrightarrow\exists
z(yz=x))).$\\

$({\bf N})_{n}$: normality up to degree $n\in\mathbb{N}^{*}$
(i.e., being domain and for every $m\in\mathbb{N}^{*}$, $m\leq n$,
$\forall x, y, z_1, \ldots, z_{m}$

\hspace{.8in}$(y\neq 0\wedge
x^{m}+z_{1}x^{m-1}y+\ldots+z_{m-1}xy^{m-1}+z_{m}y^{m}=0\longrightarrow\exists
z(yz=x))).$

It is clear that any domain satisfies $(N)_{1}$.\\

$\bf{GCD}$: having greatest common divisor (i.e., the usual axioms
for being a domain plus
\begin{center}
$\forall x,y (x=y=0 \vee\exists z (z|x \wedge z|y \wedge (\forall
t ((t|x \wedge t|y)\rightarrow t|z))))$,
\end{center}
where $x|y$ is an abbreviation for $\exists t (t\cdot x=y)$).\\

$\bf{Bez}$: the usual axioms for being a domain plus the {\em
Bezout property}:

\begin{center}
$\forall x,y \exists z,t ((xz+yt)|x \wedge (xz+yt)|y)$,
\end{center}
namely, every finitely generated ideal is principal.\\

It is known that Bez $\vdash$ GCD $\vdash$ N, and OI $\nvdash$
OI + N $\nvdash$ OI + GCD $\nvdash$ OI + Bez (Smith \cite{smith2},
 Lemmas 1.9 and 1.10).\\

Also we will need another algebraic property, though it is not
first-order expressible:\\

$\bf{DCC}$ : let $M$ be a domain. $M$ has the {\em divisor chain
condition} (DCC) if $M$ contains no infinite sequence of elements
$a_{0}, a_{1}, a_{2},\ldots$ such that each $a_{i+1}$ is a proper
divisor of $a_{i}$ (i.e., $a_{i}/a_{i+1}$ is a nonunit).\\

Let $M$ be an ordered domain (resp. a domain), then $RC(M)$ (resp.
$AC(M)$) will denote the real closure (resp. the algebraic closure)
of its fraction field. It is well known
that $AC(M)=RC(M)[\sqrt{-1}]$. Let $p\in\mathbb{N}^{*}$ and $F$ be an
ordered field (resp. a field), we define the $p$-real closure (resp.
the $p$-algebraic closure) of $F$, denoted by $RC_{p}(F)$ (resp.
$AC_{p}(F)$), to be the smallest subfield of $RC(F)$ (resp. $AC(F)$)
containing $F$ such that every polynomial of degree $\leq p$ with
coefficients in $RC_{p}(M)$ (resp. $AC_{p}(F)$) which has a root in
$RC(F)$ (resp. $AC(F)$) also has a root in $RC_{p}(F)$ (resp.
$AC_{p}(F)$). Similarly if $M$ be an ordered domain (resp. a
domain), then $RC_{p}(M)$ (resp. $AC_{p}(M)$) will denote the
$p$-real closure (resp. the $p$-algebraic closure) of its fraction
field. It can be shown that $AC_{p}(M)=RC_{p}(M)[\sqrt{-1}]$.
Similar to real closed fields and algebraic closed fields, it
is also easily seen that:\\

(1) {\em If $P(x)$ is a polynomial of degree $\leq p$ with the
coefficients in $RC_{p}(F)$ and $P(a)<0<P(b)$, for some $a<b$ in
$RC_{p}(F)$, then there exists a $c\in RC_{p}(F)$, such that
$a<c<b$ and $P(c)=0$.}\\

(2) {\em If $P(x)$ is a polynomial of degree $\leq p$ with the
coefficients in $AC_{p}(F)$, then $P(x)$  can be represented as a
product of linear factors with coefficients in $AC_{p}(F)$.}\\

Properties (1) and (2) can define and axiomatize the notions of
{\em $p$-real closed field} and {\em $p$-algebraic closed field},
denoted by $($RCF$)_{p}$ and $($ACF$)_{p}$, respectively.

Given two ordered domains $I\subset{K}$ we say that $I$ is an {\em
integer part} of $K$ if $I$ is discrete and for every element
$\alpha\in K$, there exists an element $a\in I$ such that
$0\leq\alpha-a< 1$. We call $a$, the {\em integer part} of $\alpha$,
and sometimes denote it by $[\alpha]_{I}$. Shepherdson and
Boughattas characterized models of (OI)$_{p}$, in terms of $p$-real
closed fields $(1\leq p\leq\omega)$:

\begin{theorem} [Shepherdson \cite{shep}]Let $M$ be an ordered
domain. M is a model of OI iff M is an integer part of RC(M).
\end{theorem}
\begin{theorem}[Boughattas \cite{b,Bb}]Let $M$ be an ordered domain. M is a model of
$(OI)_{p}$ iff M is an integer part of $RC_{p}(M)$.

\end{theorem}

We also need a fact from Puisseux series:

\begin{definition} Let K be a field. The following is the field of Puisseux series
in descending powers of x with coefficients in K $:$

\begin{center}
$K((x^{1/\mathbb{N}}))= \{\displaystyle \sum_{k\leq
{m}}a_{k}x^{k/r}: m\in\mathbb{Z}, r\in\mathbb{N^\ast},
a_{k}\in{K}\}.$
\end{center}
\end{definition}

\begin{theorem}[Boughattas \cite{b}] $(1\leq p\leq\omega)$ K is a $p$-real
$($resp. $p$-algebraically$)$ closed field iff
$K((x^{1/\mathbb{N}}))$ is a $p$-real $($resp. $p$-algebraically$)$
closed field.
\end{theorem}

\section{The main results}
\subsection{Proof of Theorem A}

Suppose $S$ is a subset of the set of prime numbers $\mathbb{P}$. We
present here a {\em relative to S} version of some theorems of
(Smith \cite{smith}) that is needed for proving theorem A.
Interestingly, all proofs of (Smith \cite{smith}) remain valid, if
we make routine changes which will be explained. We mention that
when $S=\mathbb{P}$, we get the original definitions and theorems.
We first define $\widehat{\mathbb{Z}}_{S}=\prod_{p\in
S}\mathbb{Z}_{p}$, where $\mathbb{Z}_{p}$ is the ring of $p$-adic
integers, and $\langle S\rangle=\{p_{1}^{\alpha_{1}}\cdots
p_{n}^{\alpha_{n}}; n\in\mathbb{N}^{*},\alpha_{i}\in\mathbb{N}$ and
$p_{i}\in S\}$. It is clear that there is the canonical embedding of
$\langle S\rangle$ in $\widehat{\mathbb{Z}}_{S}$.

Let $M$ be a model of (ZR)$_{S}$, by relativizing to $S$, we get a
(unique) $S-remainder$ homomorphism Rem :
$M\longrightarrow\widehat{\mathbb{Z}}_{S}$ given by the projective
limit of the canonical homomorphism
\begin{center}
$\psi_{n}:M\longrightarrow M/nM\cong\mathbb{Z}/n\mathbb{Z}$
\end{center}
for $n\in\langle S\rangle$. See (Macintyre-Marker \cite{MM}, Lemma
1.3).

Now we give the $S$-relativization of the so called
$\widehat{\mathbb{Z}}$-construction. Let $M$ be a discretely
ordered ring with
$\varphi:M\longrightarrow\widehat{\mathbb{Z}}_{S}$ a homomorphism
and assume that all standard primes remain prime in $M$. We form a
new ring $M_{\varphi,S}=\{a/n;a\in M,n\in\langle S\rangle$ and
$n|\varphi(a)$ in $\widehat{\mathbb{Z}}_{S}\}$. We extend
$\varphi$ to $M_{\varphi,S}$ in the obvious way. We say that
$M_{\varphi,S}$ is obtained from $M$ by the
$\widehat{\mathbb{Z}}_{S}-construction$. By relativizing the proof
of (Macintyre-Marker \cite{MM}, Lemma 3.1) we get:

\begin{lemma}
$M_{\varphi,S}\models (ZR)_{S}$.
\end{lemma}

Parsimony of homomorphisms plays a very important role in Smith's
constructions. Therefore we have the following definition:

\begin{definition}
Let $M$ be a discretely ordered ring with
$\varphi:M\longrightarrow\widehat{\mathbb{Z}}_{S}$ a homomorphism,
where $\varphi$ is the projective limit of the homomorphism
$\psi_{n}:M\longrightarrow\mathbb{Z}/n\mathbb{Z}$ for $n\in\langle
S\rangle$. We say that $\varphi$ is $S$-\textsl{parsimonious} if
for each nonzero $a\in M$ there are only finitely many
$n\in\langle S\rangle$ such that $\psi_{n}(a)=0$.
\end{definition}

The following lemma asserts that the
$\widehat{\mathbb{Z}}_{S}$-construction preserves parsimony.

\begin{lemma}
 If $\varphi:M\longrightarrow\widehat{\mathbb{Z}}_{S}$
is S-parsimonious, then the extension of $\varphi$ to
$M_{\varphi,S}$ is S-parsimonious.
\end{lemma}

\begin{proof}
The proof is the $S$-relativization of Smith's proof of Lemma 5.1.
in \cite{smith}. Let $0\neq a/n\in M_{\varphi,S}$, where $a\in M$,
$n\in\langle S\rangle$. Suppose $\psi_{k}(a/n)=0$, for a
$k\in\langle S\rangle$. Since $M_{\varphi,S}$ is a model of
(ZR)$_{S}$, we have $k|a/n$ in $M_{\varphi,S}$, so in particular
$k|a$ in $M_{\varphi,S}$. Thus $\psi_{k}(a)=0$. Since
$\varphi:M\longrightarrow\widehat{\mathbb{Z}}_{S}$ is
$S$-parsimonious, there are only finitely many possibilities for
$k\in\langle S\rangle$.
\end{proof}

The following theorem says that in the presence of having a
$S$-parsimonious map the $\widehat{\mathbb{Z}}_{S}$-construction
preserves GCD and DCC.

\begin{theorem}
Let M be a discretely ordered ring with the GCD (DCC). Let
$\varphi:M\longrightarrow\widehat{\mathbb{Z}}_{S}$ be
S-parsimonious and in the DCC case the standard primes remain
prime in M. Then $M_{\varphi,S}$ has the GCD (DCC).
\end{theorem}
\begin{proof}
We leave the proof to the reader as an easy and instructive
exercise to adopt Smith's proofs of Theorems 5.3. and 5.5. in
\cite{smith}. Just replace everywhere in the proof,
$\mathbb{Z}$-ring by a model of (ZR)$_{S}$,
$\varphi:M\longrightarrow\widehat{\mathbb{Z}}$ by
$\varphi:M\longrightarrow\widehat{\mathbb{Z}}_{S}$, parsimonious
by $S$-parsimonious, $M_{\varphi}$ by $M_{\varphi,S}$, and check
that the arguments remain valid!
\end{proof}

Transcendental extensions preserve GCD and DCC.

\begin{theorem}[Smith \cite{smith}, Theorems 6.8. and 6.10.] Let M be a
GCD (DCC) domain and suppose x is transcendental over M. Then
$M[x]$ is a GCD (DCC) domain.
\end{theorem}

By the same adaptation of Theorem 6.12. of (Smith \cite{smith}),
we see that $S$-parsimonious maps can be extended to
transcendental extensions. More precisely:

\begin{theorem}Let M be
a countable model of $(ZR)_{S}$ and suppose the remainder
homomorphisms $\varphi:M\longrightarrow\widehat{\mathbb{Z}}_{S}$
is S-parsimonious. Let $x$ be transcendental over $M$ and suppose
$M[x]$ is discretely ordered $($and this ordering restricts to the
original ordering on M $)$. Then $\varphi$ can be extended to a
S-parsimonious $\varphi:
M[x]\longrightarrow\widehat{\mathbb{Z}}_{S}$, such that
$\varphi(x)$ is a unit of $\widehat{\mathbb{Z}}_{S}$.
\end{theorem}

We will need in this paper to consider the property of {\em
factoriality} (a factorial domain has the property that any
nonunit has a factorization into irreducible elements, and this
factorization is unique up to units). We will use the following
theorem:

\begin{theorem}[Smith \cite{smith}, Theorem 1.5.]
M is factorial iff M has both of the GCD property and DCC.
\end{theorem}

In order to gain a Bezout domain the F-construction in
Macintyre-Marker paper \cite{MM} has a crucial role. By combining
Theorems 8.5 and 8.7 from (Smith \cite{smith}), Lemma 3.26 of
(Macintyre-Marker \cite{MM}) and its proof, we have:

\begin{theorem}
Let M be a discretely ordered domain with DCC (GCD) and suppose
$v,w\in M$ are primes and x is larger than any element of M. Let
$M^{*}=M[x,\frac{1-xv}{w}]$. Then $M^{*}$ is a discretely ordered
domain with DCC (GCD).
\end{theorem}

In the following theorem we see that $S$-parsimony can be extended
in F-constructions:

\begin{theorem}
Let $M$ be a  countable model of $(ZR)_{S}$ and the remainder
homomorphism $\varphi:M\longrightarrow\widehat{\mathbb{Z}}_{S}$ is
$S$-parsimonious. Let $v,w\in M$ be primes of $M$ and $w$ is
nonstandard. Suppose $x$ is transcendental over $M$, and the
discrete ordering of M extends to discrete ordering on
$M^{*}=M[x,\frac{1-xv}{w}]$. Then $\varphi$ can be extended to
$S$-parsimonious
$\varphi:M^{*}\longrightarrow\widehat{\mathbb{Z}}_{S}$, such that
$\varphi(x)$ is a unit of $\widehat{\mathbb{Z}}_{S}$.
\end{theorem}

\begin{proof}
See the proof of Theorem 8.9. of (Smith \cite{smith}).
\end{proof}

The next theorem guarantees the preservation of the GCD property
and DCC in chains constructed by alternative applications of the
F-construction and the $\widehat{\mathbb{Z}}_{S}$-construction via
parsimonious maps. We express the theorems in a more restricted
and more suitable form which is adequate for us:

\begin{theorem}
Suppose $M_{0}$ is a (GCD) DCC countable model of $(ZR)_{S}$ and
there is $S$-parsimonious
$\varphi:M_{0}\longrightarrow\widehat{\mathbb{Z}}_{S}$. Let
$\{M_{i} : i\in{\mathbb{N}}\}$ be a chain of discretely ordered
domains such that $M_{2i+1}$ is constructed from $M_{2i}$ by the
$\widehat{\mathbb{Z}}_{S}$-construction, and $M_{2i+2}$ is
constructed from $M_{2i+1}$ by the F-construction. In addition we
suppose that in the DCC case, in the whole process of extending
rings at most finitely many irreducibles have been killed (this means that
only finitely many irreducibles
will become reducible in later stages). Then
$M=\bigcup_{i\in{\mathbb{N}}}{M_{i}}$ is a model of (GCD) DCC.
\end{theorem}

\begin{proof}
See Theorems 9.4. and 9.8. in (Smith \cite{smith}).
\end{proof}

By the following series of easy lemmas, we will not worry about
DCC in our chain of models in the proof of Theorem A:

\begin{lemma}[Smith \cite{smith}, Lemma 3.8.] Let M be a GCD domain.
Then $p\in M$ is irreducible iff it is prime.
\end{lemma}

Of course the following lemma needs an easy $S$-adaptation of
Lemma 3.2 in (Macintyre-Marker \cite{MM}):

\begin{lemma}
Let $M\models DOR$ and
$\varphi:M\longrightarrow\widehat{\mathbb{Z}}_{S}$ be a ring
homomorphism and assume that all standard primes remain prime in
M. If $q\in M$ is irreducible and $\varphi(q)$ is unit in
$\widehat{\mathbb{Z}}_{S}$, then q is irreducible in
$M_{\varphi,S}$
\end{lemma}

\begin{lemma}[Macintyre-Marker \cite{MM}, Lemma 3.27]If q is
irreducible in M, then q is irreducible in $M^{*}$, constructed in
Theorem 3.8 (by the F-construction).
\end{lemma}

Now we have gathered all preliminaries to prove Theorem A:\\

\noindent {\bf Theorem A.} {\em Suppose S is a nonempty subset of
$\mathbb{P}$ and $q$ is prime number such that $q\notin S$, then
$(ZR)_{S}$ + $Bez$ $\nvdash$
 $(ZR)_{q}$.}\\

\begin{proof} We do a suitable and modified version of
Smith's process to construct a Bezout model of open induction
(Smith \cite{smith} Theorem 10.7.). We shall inductively construct
an $\omega$-chain of models $M_{i}$ such that
$\bigcup_{i}M_{i}=M_{\omega}$ will be a model of $(ZR)_{S}$ +
$Bez$+ $\neg(ZR)_{q}$. We work inside the ordered field
$\mathbb{Q}(x_{1},...,x_{i},...)$ so that for each $i\in\omega$,
$x_{i+1}$ is larger than any element of
$\mathbb{Q}(x_{1},...,x_{i})$ and $x_{1}$ is infinitely large. We
will do the F-construction at odd stages and
the $\widehat{\mathbb{Z}}_{S}$-construction at even stages.

Take $M_{0}=\mathbb{Z}$ together the natural remainder
$S$-parsimonious homomorphism
$\varphi:M_{0}\longrightarrow\widehat{\mathbb{Z}}_{S}$. Let us show
what we do at stages $2k+1$. Suppose $M_{2k}$ and a $S$-parsimonious
map $\varphi:M_{2k}\longrightarrow\widehat{\mathbb{Z}}_{S}$, have
been constructed. At this stage we consider a pair of distinct
primes $v$ and $w$ belonging to $M_{2k}$ such that $w$ is
nonstandard. (Of course we do this in such a way that every such
pair of primes in $M_{\omega}$ will have been considered at some
stage $2k+1$). Thus $(v,w)=1$ in $M_{2k}$. We define
$M_{2k+1}=M_{2k}[x_{k},\frac{1-x_{k}v}{w}]$ according to Theorem
3.8. Suppose $y_{k}=\frac{1-x_{k}v}{w}$, then we have
$x_{k}v+y_{k}w=1$ in $M_{2k+1}$. So $(v,w)_{B}=1$ in $M_{2k+1}$.
($(v,w)_{B}$ is the {\em Bezout greatest common divisor} of $v$ and $w$, it means that $(v,w)_{B}|v$ and  $(v,w)_{B}|w$ and there exist $r$ and $s$ in $M_{2k+1}$ such that $rv+su=1$).
We refer to (Smith \cite{smith}, Section 3) for the basic related
definitions and theorems. By Theorem 3.9 $\varphi$ is extended to
a $S$-parsimonious map $\varphi:
M_{2k+1}\longrightarrow\widehat{\mathbb{Z}}_{S}$. At stage $2k+2$,
we employ Lemma 3.1 and define $M_{2k+2}=(M_{2k+1})_{\varphi,S}$
which is a model of $(ZR)_{S}$. Lemma 3.3 gives us the desired
parsimonious extensions $\varphi:
M_{2k+2}\longrightarrow\widehat{\mathbb{Z}}_{S}$. Since $(ZR)_{S}$
is a $\forall\exists$-theory, then it is preserved in chains,
therefore $M_{\omega}$$\models$$(ZR)_{S}$.

Now we show that $M_{\omega}$ is a Bezout domain. The proof is
similar to (Smith \cite{smith}, Theorem 10.7) with a minor change.
By Theorems  3.4 and 3.8, each $M_{i}$ has the GCD and DCC, so by
Theorem 3.10 $M_{\omega}$ has both the GCD and DCC (by Lemmas 3.11,
3.12 and 3.13 we know that no irreducible is killed) and from
Theorem 3.7 we conclude that $M_{\omega}$ is a factorial domain. In
order to show that $M_{\omega}$ is a Bezout domain, by considering
the fact that $M_{\omega}$ has the GCD property, it suffices to
prove that any two elements of $M_{\omega}$ has the Bezout greatest
common divisor. Let $a,b\in M_{\omega}$ and let $c=(a,b)$ in
$M_{\omega}$. We can assume $a,b>1$. Let $a=a^{'}c$, $b=b^{'}c$ in
$M_{\omega}$. So $(a^{'},b^{'})=1$ in $M_{\omega}$. Since
$M_{\omega}$ is factorial, we can write
$a^{'}=m{P_{1}}^{e_{1}}\ldots{P_{k}}^{e_{k}}$ and
$b^{'}=n{Q_{1}}^{f_{1}}\ldots{Q_{l}}^{f_{l}}$, where
$m,n\in\mathbb{N}$ are nonzero, $k,l\geq 0$ and the $P_{i},Q_{j}$
are nonstandard primes such that $P_{i}\neq Q_{j}$ for all $i,j$. We
will show that $(a^{'},b^{'})_{B}=1$. Clearly $(m,n)_{B}=1$. Suppose
$m={g_{1}}^{v_{1}}\ldots{g_{r}}^{v_{r}}$ and
$n={h_{1}}^{w_{1}}\ldots{h_{s}}^{w_{s}}$ are the prime
factorizations of $m,n$ in $\mathbb{N}$. By the F-construction every
one of $(P_{i},g_{j})_{B}=1$, $(Q_{i},h_{j})_{B}=1$ and
$(P_{i},Q_{j})_{B}=1$, occur at some odd stage of our construction.
Therefore by iterated applications of (Smith \cite{smith}, Lemma
3.4), we conclude that $(a^{'},b^{'})_{B}=1$. By (Smith
\cite{smith}, Lemma 3.4), we have $c=(a,b)_{B}$ at some odd stage
and then using (Smith \cite{smith}, Lemma 3.7) we ensure that
$c=(a,b)_{B}$ in $M_{\omega}$. This completes the proof of the
Bezoutness of $M_{\omega}$.

Note that in the original proof of Smith (\cite{smith}, Theorem
10.7) he just considers pairs of nonstandard primes and doesn't need
to consider pairs of primes such that one is standard and the other is
nonstandard. Since his chain of domains are ZR-rings, this gives
automatically the Bezout greatest common divisor for such pairs. But
as we want ZR to fail in our model, we are forced to consider pairs
of standard and nonstandard primes in the F-construction, as well.

Now we show that $(ZR)_{q}$ fails in $M_{\omega}$. We first observe
that in the first step of our construction, namely, when passing
from $M_{0}=\mathbb{Z}$ to $M_{1}$, there is no nonstandard prime in
$M_{0}$. So $M_{1}$ is just $\mathbb{Z}[x_{1}]$ and we have no
$y_{1}$. On the other hand from the construction it is evident that
elements of $M_{\omega}$ are of the form
$f(x_{1},x_{2},y_{2},...,x_{k},y_{k})$, for some $k$, where $f$
is a polynomial with the coefficients in the set
$\mathbb{Z}_{\langle S\rangle}=\{a/k;a\in\mathbb{Z}$ and
$k\in\langle S\rangle\}$. Now for a contradiction, suppose $M_{\omega}$
is a model of $(ZR)_{q}$. Then there is a $b\in M_{\omega}$ such
that $x_{1}=bq+r$ with $0\leq r< q$. Take
$b=f(x_{1},x_{2},y_{2},...,x_{k},y_{k})$, so we have
$x_{1}=f(x_{1},x_{2},y_{2},...,x_{k},y_{k})q+r$. Observe that
$x_{2}$, $y_{2}$, $\ldots$, $x_{k}$, $y_{k}$ are transcendental over
$\mathbb{Q}(x_{1})$, then $f$ does not depend on them, so we can
assume $x_{1}=f(x_{1})q+r$. Since $x_{1}$ is also transcendental
over $\mathbb{Q}$, it follow that the degree of $f$ must be one.
Thus $f(x_{1})=ax_{1}$ and $a\in\mathbb{Z}_{\langle S\rangle}$. So
$x_{1}=ax_{1}q+r$ and then $x_{1}(1-aq)=r$, which implies that
$a=1/q$ and this is in contradiction with $a\in\mathbb{Z}_{\langle
S\rangle}$, since $q\notin\langle S\rangle$.
 \end{proof}

 \subsection{Proof of Theorem B}

Now we prove:\\

{\bf Theorem B.} {\em Suppose S is a nonempty subset of
$\mathbb{P}$ and $q$ is a prime number such that $q\notin S$, then
$(ZR)_{S} + (N)_{n} + \neg (ZR)_{q} + \neg (N)_{n+1}$ is
consistent.}\\

\begin{proof} In \cite{moh2} we proved that if $n\neq 3$, there
is a $\lambda$ which is real algebraic of degree $n+1$ over
$\mathbb{Q}$ and doesn't belong to $RC_{n}(\mathbb Q)$. Now
suppose $x$ is an infinitely large element. For $n\neq3$, fix $\lambda$ as
above. For $n=3$ we choose $\lambda$ as a root of an irreducible
polynomial of degree 4 such that $\lambda\notin RC_{2}(\mathbb
Q)$. Let $A$ be the ring of integers of the algebraic number field
$\mathbb{Q}(\lambda)$. Form $A_{\langle S\rangle}=\{a/k;a\in A$
and $k\in\langle S\rangle\}$. It is an elementary fact from algebraic
number theory that $A$ is a normal ring. Since $A_{\langle
S\rangle}$ is a localization of $A$ relative to a multiplicative
set, then it is also normal. Let $M=\mathbb{Z}[rx;r\in A_{\langle
S\rangle}]$. We claim that $M$ witnesses Theorem B. It is obvious
that $M\vDash(ZR)_{S}$. By an argument similar to the last
paragraph of the proof of theorem A, it is easily shown that
$M\vDash\neg (ZR)_{q}$.

Now we prove $M\vDash\neg(N)_{n+1}$. Let $v\in\mathbb{N}$ be such
that $v\lambda$ is an algebraic integer. Suppose
$P(t)\in\mathbb{Z}$ is its minimal polynomial of degree $n+1$
which is monic. Obviously $v\lambda x\in M$. But we have
$P(v\lambda x/x)=0$, while $v\lambda\notin M$. So $M$ is not a
model of $(N)_{n+1}$.

It remains to show that $M\vDash(N)_{n}$. Let $u,v$ be nonzero
elements of $M$ such that
\begin{center}
$(u/v)^{s}+z_{1}(u/v)^{s-1}+\cdots+z_{s}=0$\,\,\,\,\,\,\,$(z_{1},\ldots,z_{s}\in
M, s\leqslant n).$
\end{center}
We will show that $u/v\in M$. Notice that elements of $M$ are those
elements of $A_{\langle S\rangle}[x]$ with integer constant
coefficient. $ A_{\langle S\rangle}$ is normal, so is $ A_{\langle
S\rangle}[x]$. Thus $u/v\in A_{\langle S\rangle}[x]$. On the other
hand, since $\mathbb{Q}(\lambda)[x]$ is a factorial ring, $u/v$ can be
written as:
\begin{center}

$u/v=\rho\prod_{i\in I}P_{i}\prod_{j\in J}Q_{j}$,
\end{center}
in which $\rho\in\mathbb{Q}(\lambda)$, the $P_{i}$'s are irreducible
in $\mathbb{Q}(\lambda)[x]$, without constant coefficient and
$Q_{j}$'s are irreducible in $\mathbb{Q}(\lambda)[x]$ with the
constant coefficient one. If $I$ is nonempty, then $\rho\prod_{i\in
I}P_{i}\prod_{j\in J}Q_{j}$ has no constant coefficient and thus
$u/v\in M$. Now suppose $I=\O$. Put $x=0$ in $u$, $v$, $z_{1}$,
$\ldots$, $z_{s}$. Therefore $\rho$ is an algebraic integer with the
degree, equal or less than $n$ over $\mathbb{Z}$. We show it is one.
If $n=1$ there is nothing to prove. If not, we have
$[\mathbb{Q}(\lambda): \mathbb{Q}(\rho)]< n+1$. But

\begin{center}

$[\mathbb{Q}(\lambda): \mathbb{Q}(\rho)][\mathbb{Q}(\rho):
\mathbb{Q}]=[\mathbb{Q}(\lambda): \mathbb{Q}]=n+1$.

\end{center}
Then $[\mathbb{Q}(\lambda): \mathbb{Q}(\rho)]$ divides
$[\mathbb{Q}(\lambda): \mathbb{Q}]=n+1$. So we have a chain of
field extensions,
$\mathbb{Q}\subset\mathbb{Q}(\rho)\subset\mathbb{Q}(\lambda)$ such
that $[\mathbb{Q}(\lambda): \mathbb{Q}(\rho)]\leq n-1$ and
$[\mathbb{Q}(\rho): \mathbb{Q}]\leq n-1$. This implies that
$\lambda\in RC_{n-1}(\mathbb Q)$ which is in contradiction with
the choice of $\lambda$. Hence $\rho$ is an algebraic integer of
degree one. So $\rho\in\mathbb{Z}$ and this implies that $u/v\in
M$, which means that $M$ is model of $(N)_{n}$. This completes the
proof of Theorem B.
\end{proof}

\subsection{Proofs of Theorems C and D.}
In order to demonstrate Theorem C, we need a generalization of a
theorem of Boughattas. In (\cite{Bb}, Theorem V.1.) Boughattas proved
that every saturated ordered field admits a normal integer part. But
we show that:

\begin{lemma}Every $\omega_{1}$-saturated ordered field admits a
Bezout integer part.
\end{lemma}

\begin{proof}
({\em Sketch}) Suppose $K$ is an $\omega_{1}$-saturated ordered
field. Boughattas \cite{Bb} in a series of three Lemmas: {\em
Pricipal}, {\em Integer Part} and {\em Construction}, showed that
we can build an $\omega_{1}$-chain of countable discretely ordered
rings $M_{i}, i<\omega_{1}$ such that
$M=\bigcup_{i<\omega_{1}}M_{i}$ is an integer part of $K$.
Furthermore he considers an arbitrary subset $\Lambda\subset K$ of
real algebraic elements which plays a role in the construction of the
$M_{i}$'s. Varying $\Lambda$ gives us various kinds
of integer parts. When $\Lambda=\O$, we obtain a normal integer
part and it is implicit in the paper that in this case, the $M_{i}$'s
are obtained by alternative applications of the
Wilkie-construction and the $\widehat{\mathbb{Z}}$-construction.
But it must be noticed that even in this case the procedure of
doing the $\widehat{\mathbb{Z}}$-construction is different from
the original one, because it is no longer assumed that the ground field is dense in
its real closure. To gain a Bezout integer part, we observe that
we can do the procedure of the Theorems 10.7 and 10.8 of Smith
\cite{smith} inside $K$. In this procedure we need the extra
F-construction. Since any $M_{i}, i<\omega_{1}$ is countable and
$K$ is $\omega_{1}$-saturated, then there is always an element
$b_{i}$ in $K$ which is larger than any element of $M_{i}$. By Lemma
3.26 of (Macintyre-Marker \cite{MM}) we are sure that we can do
the F-construction.  To obtain an integer part of $K$, suppose
$(b_{\alpha},\alpha<\omega_{1})$ be an enumeration of elements of
$K$. Let $M_{i}$ has been constructed and at step $i+1$ we want to
do the Wilkie-construction. We seek the least ordinal
$\alpha_{i}$, such that $b_{\alpha_{i}}$ has not an integer part
in $M_{i}$. Then by combining the Integer Part Lemma of Boughattas
\cite{Bb} with the $\widehat{\mathbb{Z}}$-construction, we obtain
$M_{i+1}$ with its parsimonious homomorphism extension to
$\widehat{\mathbb{Z}}$, such that $b_{\alpha_{i}}$ has an integer
part in $M_{i+1}$. Also suppose $M_{j}$ has been constructed and at
stage $j+1$ we want to do the F-construction. We seek the
least ordinal $\alpha_{j}$, such that $b_{\alpha_{j}}$ is larger
than any element of $M_{j}$. Then the F-construction can be done
at this step. At limit stages we take union. Moreover, Lemma 9.1,
Theorem 9.4 and Theorem 9.8 of (Smith \cite{smith}) will guarantee
preserving parsimony of homomorphisms and factoriality at limit
stages of length $\leq\omega_{1}$. Now there is no obstacle for
$M=\bigcup_{i<\omega_{1}}M_{i}$ to be a Bezout integer part of
$K$.
\end{proof}

\noindent{\bf Theorem C.} $Bez$ + $(OI)_{p}$ $\nvdash$
$(OI)_{p+1}$, {\em when} $p\neq 3$.\\

\noindent{\bf Proof of Theorem C.} In \cite{moh2}, we showed that if
$p\neq3$, there is an irreducible polynomial $P(t)$ of degree $p+1$
over $\mathbb{Q}$ such that $P(t)$ has no root in
$RC_{p}(\mathbb{Q})$. For $p\geq4$, $P(t)$ was a polynomial with
Galois group $A_{p+1}$. It is well known that we can take $P(t)$ as
a monic polynomial with integer coefficients such that $P(0)<0$. Let
$T$ be the following theory in the language of ordered field with
the additional constant symbol $a$:
\begin{center}
$T\equiv(RCF)_{p}+ \{a>k; k\in\mathbb{N}\}+ \forall y
\neg(Q(y)\leq0<Q(y+1))$,
\end{center}
where $Q(y)=a^{p+1}P(y/a)$. We show that the field of Puisseux power
series $RC_{p}(\mathbb{Q})((x^{1/\mathbb{N}}))$ is a model of $T$,
when interpreting $a$ by $x$. Clearly by Theorem 2.4,
$RC_{p}(\mathbb{Q})((x^{1/\mathbb{N}}))$ is a $p$-real closed field.
Also on the contrary suppose that there exists $y$ $\in$
$RC_{p}(\mathbb{Q})((x^{1/\mathbb{N}}))$ such that
$(Q(y)\leq0<Q(y+1))$. Therefore $P(y/x)\leq0<P((y+1)/x)$. It is
easily seen that deg$_{x}(y/x)$ must be zero. So let
$y/x=\lambda+\sum_{-\infty<i<0}c_{i}x^{i/q}$ in
$RC_{p}(\mathbb{Q})((x^{1/\mathbb{N}}))$. This leads to
$P(\lambda)=0$, but $\lambda\in RC_{p}(\mathbb{Q})$, which is in
contradiction to the choice of $P(t)$. So
$RC_{p}(\mathbb{Q})((x^{1/\mathbb{N}}))\vDash T$.

Now that $T$ is consistent, let $K$ be an $\omega_{1}$-saturated
model of $T$. By Lemma 3.14 $K$ has a Bezout integer part. Call it
$M$. Since $K\vDash(RCF)_{p}$, then $M\vDash(OI)_{p}$. On the other
hand there is $n\in\mathbb{N}$ such that $M\vDash
Q(0)<0<Q(n[a]_{M})$, where $[a]_{M}$ is the integer part of $a$ in
$M$. But $K\vDash\forall y \neg(Q(y)\leq0<Q(y+1))$, then
$M\vDash\forall y \neg(Q(y)\leq0<Q(y+1))$, so by (Boughattas
\cite{b}, Proposition A.I), $M\vDash\neg(OI)_{p+1}$. This ends the
proof of $Bez$ + $(OI)_{p}$ $\nvdash$ $(OI)_{p+1}$.
$\blacksquare$

Proof of Theorem D goes the same way with the exception that we must
replace Lemma 3.14 by the following Construction Lemma of
Boughattas:

\begin{theorem}[Boughattas \cite{Bb}] Suppose $K$ is a saturated
ordered field. Let $\Lambda$ be an arbitrary subset of real
algebraic elements in $K$. Then there exists $X\subset K$ such
that $X$ is algebraic independent and
$\mathbb{Z}[\{rx;r\in\mathbb{Q}[\Lambda]$ $and$ $x\in X\}]$ is an
integer part of $K$.
\end{theorem}

\noindent{\bf Theorem D.} $(OI)_{p}$ + $\neg(OI)_{p+1}$ +
$(N)_{n}$ + $\neg(N)_{n+1}$ {\em is consistent, when} $p\neq 3$.\\

\noindent{\bf Proof of Theorem D.} We work with the same theory $T$
and its saturated model $K$ as in the proof of Theorem C. Choose
$\Lambda=\{\lambda\}$ and fix $\lambda$ as in the proof of Theorem
B, namely, if $n\neq3$, $\lambda\in RC_{n+1}(\mathbb Q)$,
$\lambda\notin RC_{n}(\mathbb Q)$ and if $n=3$ choose $\lambda$ as a
root of an irreducible polynomial of degree 4 such that
$\lambda\notin RC_{2}(\mathbb Q)$. Then by Theorem 3.15, there exists $X\subset K$ such
that $K$ has the
integer part $M=\mathbb{Z}[\{rx;r\in\mathbb{Q}(\lambda)$ and $x\in
X\}]$. To show $M \models (N)_{n} + \neg(N)_{n+1}$, we can repeat the
proof of Theorem B, just replace $x$ by $X$ and replace $A_{\langle
S\rangle}$ by $\mathbb{Q}(\lambda)$. $\mathbb{Q}(\lambda)[X]$
remains factorial and normal, so the proof works. By the last
paragraph of the proof of the Theorem C, it is obvious that
$M\models(OI)_{p}+\neg(OI)_{p+1}$ $\blacksquare$

\subsection*{Acknowledgment} I would like to thank Professor Roman
Kossak for his patience and kindness during the preparation of this
paper.

\bibliography{subsystem}
\bibliographystyle{plain}
\end{document}